\documentclass[12pt]{article}%
\usepackage{amssymb,amsmath}
\usepackage{color}%
\usepackage{amsfonts}%
\usepackage{graphicx}
\providecommand{\U}[1]{\protect\rule{.1in}{.1in}}

\newtheorem{theorem}{Theorem}
\newtheorem{proposition}{Proposition}
\newtheorem{example}{Example}
\newtheorem{lemma}{Lemma}
\newtheorem{remark}{Remark}
\newtheorem{definition}{Definition}
\newtheorem{corollary}{Corollary}
\newenvironment{proof}[1][Proof]{\noindent\textbf{#1.} }{\ \rule{0.5em}{0.5em}}

\def\calf{{\cal F}}

\newcommand\id{\operatorname{id}}
\def\RR{\mathbb{R}}
\def\<{\langle}
\def\>{\rangle}
\newcommand{\eq}{\hspace*{-.4mm}&=&\hspace*{-.4mm}}

\begin{document}

\title{An integral formula for a pair of singular distributions}
\author{
Paul Popescu\thanks{Department of Applied Mathematics, University of Craiova,
Craiova 200585, Str. Al. Cuza, No, 13, Romania. e-mail: \texttt{paul$\_$p$\_$popescu@yahoo,com}}
\ and \
Vladimir Rovenski\thanks{Department of Mathematics, University of Haifa, Mount Carmel, 31905 Haifa, Israel.
\newline e-mail: \texttt{vrovenski@univ.haifa.ac.il}}
}

\date{}
\maketitle

\begin{abstract}
The paper is devoted to differential geometry
of singular distributions (i.e., of varying dimension) on a Riemannian manifold.
Such distributions are defined as images of the tangent bundle under smooth endomorphisms.
We prove the novel divergence theorem with the divergence type operator
and deduce the Codazzi type equation for a pair of singular distributions.
Tracing our Codazzi type equation yields expression of the mixed scalar curvature
through invariants of distributions, which provides some splitting results.
Applying our divergence theorem, we get the integral formula, generalizing the known one, with the mixed scalar curvature
of a pair of transverse singular distributions.

\vskip1.5mm\noindent\textbf{Keywords}:
Riemannian metric,
singular distribution,
second fundamental form,
divergence,
integral formula,
mean curvature

\vskip1.5mm \noindent\textbf{Mathematics Subject Classifications (2010)} 53C15

\end{abstract}


\section*{Introduction}

Distributions, being subbundles of the tangent bundle $TM$ on a manifold $M$,
arise in such topics of differential geometry as vector fields, submersions, fiber bundles,
Lie groups actions, \cite{BF,GHV,rov-m}, and in theoretical physics \cite{BL,Is}.
Foliations, which are defined as partitions of a manifold $M$ into collections of submanifolds-leaves
(of the same dimensi\-on in regular case), correspond to integrable distributions.
Riemannian foliations (that is having equidistant leaves) with singularities were defined by P.\,Molino \cite{Mol},
the orbit decomposition of an isometric actions of a Lie group gives an example,~\cite{abt}.
There is some interest of geometers and engineers to singular distributions, i.e., having varying dimension, e.g. \cite{DLPR}.
 We define such distributions as images of $TM$ under smooth endomorphisms $P$.
 The~paper is devoted to differential geometry of singular distributions and foliations
 (i.e., the geometrical properties depending on structural tensors) and continues the study \cite{PP1,PP0}.
In Section~\ref{sec:01}, we deduce the Codazzi equation for a pair of transverse singular distributions.
In Section~\ref{sec:02} we prove the new divergence theorem (and its modification for open Riemannian manifolds)
 with the divergence type operator, called the $P$-divergence.
We give examples with Einstein tensors and with almost contact structure and $f$-structure.
Tracing our Codazzi equation yields expression of the mixed scalar curvature $S^P_{\rm mix}$ (see \cite{BF,rov-m} for regular distributions)
through invariants of distributions, which provides some splitting results.
Integral formulas (mainly, based on the divergence theorem) provide obstructions to existence~of foliations or compact leaves of them, see survey \cite{arw}.
In~Section~\ref{sec:03}, using our $P$-divergence theorem, we obtain the integral formula
with $S^P_{\rm mix}$ of a pair of transverse singular distributions, parameterized by self-adjoint endomorphisms;
the formula generalizes the known one in \cite{Wa01}, which has many~applications.

\section{Structural tensors of singular distributions}
\label{sec:01}

Let $M$ be a smooth $n$-dimensional manifold, $TM$ -- the tangent bundle,
${\cal X}_{M}$ -- the Lie algebra of smooth vector fields on $M$, and
$\mathrm{End}(TM)$ -- smooth endomorphisms of $TM$, i.e., linear maps on the fibers of $TM$.
Let $ C^\infty(M)$ be the algebra of smooth functions on $M$.

\begin{definition}\rm
An image ${\cal D}=\Pi(TM)$ of an endomorphism $\Pi\in\mathrm{End}(TM)$ will
be called a \emph{generalized vector subbundle} of $TM$ or a \emph{singular distribution}.
Let $\Pi({\cal X}_{M})$ be an $ C^\infty(M)$-submodule of ${\cal X}_{\cal D}$ (smooth vector fields on ${\cal D}$),
i.e., sections $Y=\Pi(X)\in{\cal X}_{{\cal D}}$,
where~$X\in{\cal X}_{M}$.
\end{definition}

Let $P_{1},P_{2}$ be endomorphisms of $TM$ such that the intersection of their images is trivial,
hence ${\rm rank}\,P_1(x)+{\rm rank}\,P_2(x)\le n$ for any $x\in M$.
For example, $P_{i}$ may be projectors from $TM$ onto transverse distributions.
Given distributions ${\cal D}_{i}=P_{i}(TM)$,
put ${\cal D}=P(TM)$ for $P=P_{1}+P_{2}\in\mathrm{End}(TM)$.
Then $P(TM)={\cal D}_{1}\oplus{\cal D}_{2}$ is the generalized vector subbundle of $TM$, but not necessarily ${\cal D}_{1}\oplus{\cal D}_{2}=TM$.
A Riemannian metric $g=\<\cdot\,,\cdot\>$ on $M$ is \textit{adapted} if
\begin{equation}\label{E-P1P2}
 P_{i}P_{j}^* = P_{i}^*P_{j} = 0\quad (i\ne j).
\end{equation}
In this case, ${\cal D}_{1}\perp_{\,g}{\cal D}_{2}$ and ${\cal D}^*_{1}\perp_{\,g}{\cal D}^*_{2}$,
where ${\cal D}^*_{i}=P^*_{i}(TM)$ are generalized vector subbundles and $P^*_{1},P^*_{2}$
are adjoint endomorphisms of $TM$, i.e.,
 $\<P_i(X),Y\> = \<X,P^*_i(Y)\>$.
Similarly, we define endomorphism $P^*=P^*_{1}+P^*_{2}$ and generalized vector subbundle ${\cal D}^*=P^*(TM)$.

\begin{remark}\label{Rem-01}\rm
By \cite{PP0}, given a Riemannian metric and a singular distribution ${\cal D}_{1}$ on $M$, there are self-adjoint endomorphisms $P_1$ and $P_2$ of $TM$ such that $P_1(TM)={\cal D}_{1}$ and $P_2(TM)={\cal D}_{2}$ are smooth orthogonal distributions, the direct sum decomposition $TM=P_1(TM)\oplus P_2(TM)$ is valid on a dense subset of $M$.
We use self-adjoint endomorphisms only in the last section.
\end{remark}

Let $\nabla: (\mathfrak{X}_{M})^2\to\mathfrak{X}_{M}$ be a linear connection, that is
(for any $X_{i}, Y\in{\cal X}_{M}$ and $f\in C^\infty(M)$)
\[
 \nabla_{f X_{1}+X_{2}}Y = f\nabla_{X_{1}}Y + \nabla_{X_{2}}Y,\quad
 \nabla_{Y}(f X_{1}+X_{2}) = f\nabla_{Y}X_{1} + Y(f)\cdot X_{1} + \nabla_{Y}X_{2}.
\]
The Levi-Civita connection of $g$ is a metric and torsionless linear connection, that is
\[
 Z\<X,Y\> = \<\nabla_Z\,X, Y\> + \<X, \nabla_Z\,Y\>,\quad
 \nabla_X\,Y - \nabla_Y\,X - [X,Y] = 0.
\]

\begin{definition}\rm
Given $P_1,P_2\in\mathrm{End}(TM)$, the \emph{structural tensors} of singular distributions ${\cal D}_{i}=P_i(TM)$
are bilinear maps
 $B_{i} :({\cal X}_{M})^2\rightarrow P_{i}({\cal X}_{M})\ (i=1,2)$,
given by
\begin{subequations}
\begin{equation}\label{E-B1B2}
 {B}_{1}(Y,X):=P^*_1\nabla_{P_{1}X} P_{2}Y,\quad
 {B}_{2}(X,Y):=P^*_2\nabla_{P_{2}Y} P_{1}X.
\end{equation}
Similarly, define \emph{structural tensors} $\hat{B}_{i}$ of dual distributions ${\cal D}^*_{i}$
and auxi\-liary tensors $\check{B}_{i}$:
\begin{eqnarray}\label{E-B1B2-dual}
 && \hat{B}_{1}(Y,X):=P_1\nabla_{P^*_{1}X} P^*_{2}Y,\quad
 \hat{B}_{2}(X,Y):=P_2\nabla_{P^*_{2}Y} P^*_{1}X,\\
\label{E-B1B2-dual2}
 && \check{B}_{1}(Y,X):=P_1\nabla_{P_{1}X} P^*_{2}Y,\quad
 \check{B}_{2}(X,Y):=P_2\nabla_{P_{2}Y} P^*_{1}X.
\end{eqnarray}
\end{subequations}
\end{definition}

For self-adjoint endomorphisms $P_{i}\ (i=1,2)$, we have $B_{i}=\hat{B}_{i}=\check{B}_{i}$
and ${\cal D}^*_{i}={\cal D}_{i}$.
In~particular, for orthoprojectors $P_{i}$ (thus, $P_{i}^2=P_{i}$) from $TM$ onto
complementary orthogonal regular distributions ${\cal D}_{i}$ and the Levi-Civita connection, the structural tensors $B_{i}$ have been defined in \cite[p. 31]{rov-m}.

\begin{definition}\rm
A pair $(P_1,P_2)$ (or a tensor $P=P_1+P_2$) is \emph{allowed for a
connection} $\nabla$~if
\begin{equation*}
 b_{j}^{(i)}=0,\ i,j\in\{1,2\}.
\end{equation*}
The bilinear forms
$b_{1}^{(i)}:({\cal X}_{M})^2\rightarrow P_{2}({\cal X}_{M})$
and their dual
$b_{2}^{(i)}:({\cal X}_{M})^2\rightarrow P_{1}({\cal X}_{M})$
are given by
\begin{align*}
 & b_{1}^{(1)}(X,Y) = P^*_{2}P_{2}\nabla_{P_{1}X}P^*_{1}Y -P^*_{2}\nabla_{P_{1}X}P_{1}P^*_{1}Y,\\
 & b_{1}^{(2)}(X,Y) =P^*_{2}P_{2}\nabla_{P_{1}X_{1}}P^*_{1}Y -P_{2}\nabla_{P^*_{1}P_{1}X}P^*_{1}Y,\\
 & b_{2}^{(1)}(X,Y) =P^*_{1}P_{1}\nabla_{P_{2}X}P^*_{2}Y -P^*_{1}\nabla_{P_{2}X}P_{2}P^*_{2}Y,\\
 & b_{2}^{(2)}(X,Y)=P^*_{1}P_{1}\nabla_{P_{2}X}P^*_{2}Y -P_{1}\nabla_{P^*_{2}P_{2}X}P^*_{2}Y.
\end{align*}
\end{definition}

\begin{example}\rm
A simple example of allowed endomorphism is $P = f\operatorname{id}$, where $P = f P_{1} + f P_{2}$, $P_{i}$ are
orthoprojectors onto complementary orthogonal distributions, $\operatorname{id}$ is the identity endomorphism of $TM$ and $f$
is a real function on $M$ such that its supporting set is dense in~$M$.
More examples of singular distributions of this type, even integrable, are given in \cite{PP1}.
\end{example}

\begin{lemma}\label{lm01}
If $(P_1,P_2)$
is allowed for a metric connection $\nabla$, then for all $X,Y\in{\cal X}_{M}$ we have
\begin{eqnarray}\label{E-lemma1}
\nonumber
 & P_{2} {B}_{2}(X,Y) = \hat{B}_{2}(X, P_2 Y)
 = \check{B}_{2}(P_1 X, Y)
 ,\\
 & P_{1} {B}_{1}(Y,X) = \hat{B}_{1}(Y, P_1 X)
 =
 \check{B}_{1}(P_2 Y, X) .
\end{eqnarray}
\end{lemma}

\proof Using $b_{2}^{(2)}=0$ and \eqref{E-P1P2}, we obtain
\begin{eqnarray*}
 0\eq \<P^*_{1}P_{1}\nabla _{P_{2}Y}P^*_{2}Z,\ X\> -\<P_{1}\nabla_{P^*_{2}P_{2}Y}P^*_{2}Z,\ X\>\\
  \eq \<\nabla_{P_{2}Y}P^*_{2}Z,\ P^*_{1}P_{1}X\> -\<\nabla_{P^*_{2}P_{2}Y}P^*_{2}Z,\ P^*_{1}X\>\\
  \eq-\<P^*_{2}Z,\ \nabla_{P_{2}Y}P^*_{1}P_{1}X\> +\<P^*_{2}Z,\ \nabla_{P^*_{2}P_{2}Y}P^*_{1}X\>\\
  \eq \<\hat{B}_{2}(X, P_2 Y) -\check{B}_{2}(P_1 X, Y),\ Z\>.
\end{eqnarray*}
Similarly, using $b_{2}^{(1)}=0$, we obtain
\begin{equation*}
 0 = \<P^*_{1}P_{1}\nabla_{P_{2}Y}P^*_{2}Z -P^*_{1}\nabla_{P_{2}Y}P_{2}P^*_{2}Z,\ X\>
   =\<
   -\check{B}_{2}(P_1 X, Y) +P_{2} {B}_{2}(X, Y)
   ,\ Z\>.
\end{equation*}
Then \eqref{E-lemma1}$_{1}$ follows.
Note that \eqref{E-lemma1}$_{2}$ is dual to \eqref{E-lemma1}$_{1}$ and follows from $b_{1}^{(1)}=b_{1}^{(2)}=0$.
\hfill$\square$

\begin{definition}\label{D-TTSSR}\rm
Define maps ${R^P},S_{i},{\cal T}_{i}:({\cal X}_{M})^{4}\rightarrow C^\infty(M)$,
$i\in\{1,2\}$ by
\begin{eqnarray*}
 &&\hskip-6mm
 {\cal T}_{1}(Y,X_{1},X_{2},Z) = \<P_{2}\nabla_{P_{1}X_{1}} {B}_{2}(X_{2},Y)
 - \check{B}_{2}(\nabla_{P_{1}X_{1}}P_{1}X_{2}, Y) - \hat{B}_{2}(X_2, \nabla_{P_{1}X_{1}}P_{2}Y) ,\,Z\>, \\
 &&\hskip-6mm
 {\cal T}_{2}(Y,X_{1},X_{2},Z) = \<P_{1}\nabla_{P_{2}Y} {B}_{1}(Z,X_{1})
  - \check{B}_{1}(\nabla_{P_{2}Y}P_{2}Z, X_{1}) -\hat{B}_{1}(Z, \nabla_{P_{2}Y}P_{1}X_1),\,X_{2}\>, \\
 &&\hskip-6mm S_{1}(Y,X_{1},X_{2},Z) = \<\hat{B}_{2}(X_2, \nabla_{P_{2}Y}P_{1}X_{1}),\,Z\> , \\
 &&\hskip-6mm  S_{2}(Y,X_{1},X_{2},Z) = \<\hat{B}_{1}(Z, \nabla_{P_{1}X_1}P_{2}Y),\,X_{2}\> , \\
 &&\hskip-6mm {R}^P(Y,X_{1},X_{2},Z) =
  \<P^*_{2}\nabla_{P_{2}Y}P_{2}\nabla_{P_{1}X_{1}}P^*_{1}X_{2}
   +P_{2}\nabla_{P_{2}Y}P^*_{1}\nabla_{P_{1}X_{1}}P_{1}X_{2} \\
 &&\hskip10mm -P^*_{2}\nabla_{P_{1}X_{1}}P_{1}\nabla_{P_{2}Y}P^*_{1}X_{2}
   -P_{2}\nabla_{P_{1}X_{1}}P^*_{2}\nabla_{P_{2}Y}P_{1}X_{2}
   -P_{2}\nabla_{P^*[P_{2}Y,\ P_{1}X_{1}]}P^*_{1}X_{2},\ Z\> .
\end{eqnarray*}
\end{definition}


\begin{proposition}
If $P=P_1+P_2$ is allowed for the Levi-Civita connection,
then
$S_i,T_i$ and $R^P$ are
tensor fields and the following Codazzi type equation holds:
\begin{equation}\label{E-01}
  S_{1}+{\cal T}_{1}+S_{2}+{\cal T}_{2} + {R^P} = 0.
\end{equation}
\end{proposition}

\proof
 Let us show that all maps in Definition~\ref{D-TTSSR} are tensor fields.
 In general, we have
\begin{align*}
 & {\cal T}_{1}(Y,fX_{1},X_{2},Z) = {\cal T}_{1}(Y,X_{1},X_{2},fZ) = f\,{\cal T}_{1}(Y,X_{1},X_{2},Z), \\
 & {\cal T}_{1}(Y,X_{1},fX_{2},Z) = f\,{\cal T}_{1}(Y,X_{1},X_{2},Z)
 + P_{1}X_{1}(f)\<\underline{P_{2} {B}_{2}(X_{2},Y)
 -\check{B}_{2}(P_{1}X_{2}, Y)}
 ,\, Z\>, \\
 & {\cal T}_{1}(fY,X_{1},X_{2},Z) = f\,{\cal T}_{1}(Y,X_{1},X_{2},Z)
 + P_{1}X_{1}(f)\<\underline{P_{2} {B}_{2}(X_{2},Y) - \hat{B}_{2}(X_2, P_2 Y)}
 ,\, Z\>
\end{align*}
for any $f\in C^\infty(M)$, using \eqref{E-lemma1}$_1$ to cancel underlying terms;
and similarly, for ${\cal T}_{2}$.
From the above and \eqref{E-lemma1} follows that ${\cal T}_{i}$ are tensors.
Using $(1,2)$-tensor fields $\tilde{B}_{i}$,
\begin{equation*}
 \tilde{B}_{1}(Y,X) := P_{1}\nabla _{X}P_{2}^{\ast}Y,\quad
 \tilde{B}_{2}(X,Y) := P_{2}\nabla _{Y}P_{1}^{\ast}X,
\end{equation*}
we obtain
\begin{equation*}
 S_{1}(Y,X_{1},X_{2},Z) = \big\langle \tilde{B}_{2}(X_{2}, {B}_{2}(X_{1},Y)) ,Z\big\rangle , \quad
 S_{2}(Y,X_{1},X_{2},Z) = \big\langle \tilde{B}_{1}(Z, {B}_{1}(Y,X_{1})) ,X_{2}\big\rangle .
\end{equation*}
Thus, $S_{i}:\mathcal{X}(M)^{4}\rightarrow \mathcal{X}(M)$ are tensor fields.
 We~have, using also that $\nabla$ is torsion-free,
\begin{align*}
 & ({\cal T}_{1} + S_{1})(Y,X_{1},X_{2},Z)=\<P_{2}\nabla_{P_{1}X_{1}} {B}_{2}(X_{2},Y)
 -P_{2}\nabla_{P_{2}Y}P^*_{1}\nabla_{P_{1}X_{1}}P_{1}X_{2} \\
 & -P_{2}\nabla_{P^*_{2}\nabla_{P_{1}X_{1}}P_{2}Y}P^*_{1}X_{2}
 +P_{2}\nabla_{P^*_{2}\nabla_{P_{2}Y}P_{1}X_{1}}P^*_{1}X_{2},\ Z\> \\
 & =\<P_{2}\nabla_{P_{1}X_{1}} {B}_{2}(X_{2},Y)
 -P_{2}\nabla_{P_{2}Y}P^*_{1}\nabla_{P_{1}X_{1}}P_{1}X_{2}
 -\hat{B}_{2}(X_2, [ P_{1}X_{1},P_{2}Y]),\ Z\>,
\end{align*}
and similarly, using \eqref{E-lemma1}$_2$ to cancel underlying terms on the last step,
\begin{align*}
 & ({\cal T}_{2}+S_{2})(Y,X_{1},X_{2},Z) \\
 & = \<\nabla_{P_{2}Y} {B}_{1}(Z, X_{1}),\ P^*_{1}X_{2}\> -\<\nabla_{P_{1}X_{1}}P^*_{2}\nabla_{P_{2}Y}P_{2}Z,\ P^*_{1}X_{2}\>
  +\<\hat{B}_{1}(Z, [{P_{1}X_1},\ P_{2}Y]),\,X_{2}\> \\
 &=P_{2}Y\< P_{1} {B}_{1}(Z, X_{1}),\ X_{2}\> -\<{B}_{1}(Z, X_{1}),\ \nabla_{P_{2}Y}P^*_{1}X_{2}\>
  +\< P^*_{2}\nabla_{P_{2}Y}P_{2}Z,\ \nabla_{P_{1}X_{1}}P^*_{1}X_{2}\> \\
 &+\< P^*_{2}Z,\ \nabla_{P^*_{1}[P_{2}Y,\ P_{1}X_{1}]}P^*_{1}X_{2}\> \\
 &=P_{2}Y\< P_{1} {B}_{1}(Z, X_{1}),\ X_{2}\>
  +\<\nabla_{P_{1}X_{1}}P_{1}\nabla_{P_{2}Y}P^*_{1}X_{2},\ P_{2}Z\>
  +P_{2}Y\< P_{2}Z,\ P_{2}\nabla_{P_{1}X_{1}}P^*_{1}X_{2}\> \\
 &-\< P_{2}Z,\ \nabla_{P_{2}Y}P_{2}\nabla_{P_{1}X_{1}}P^*_{1}X_{2}\>
  +\< P^*_{2}Z,\ \nabla_{P^*_{1}[P_{2}Y,\ P_{1}X_{1}]}P^*_{1}X_{2}\>\\
 &= P_{2}Y\<\underline{P_{1} {B}_{1}(Z, X_{1})
 - \check{B}_{1}(P_{2}Z, X)}
 ,\ X_{2}\>
  +\<\nabla_{P_{1}X_{1}}P_{1}\nabla_{P_{2}Y}P^*_{1}X_{2},\ P_{2}Z\>  \\
 &-\<\nabla_{P_{2}Y}P_{2}\nabla_{P_{1}X_{1}}P^*_{1}X_{2},\ P_{2}Z \>
  -\<\nabla_{P^*_{1}[P_{1}X_{1},P_{2}Y]}P^*_{1}X_{2},\ P^*_{2}Z \>\\
 &= \<P^*_{2}\nabla_{P_{1}X_{1}}P_{1}\nabla_{P_{2}Y}P^*_{1}X_{2}
  -P^*_{2}\nabla_{P_{2}Y}P_{2}\nabla_{P_{1}X_{1}}P^*_{1}X_{2}
  -P_{2}\nabla_{P^*_{1}[P_{1}X_{1},P_{2}Y]}P^*_{1}X_{2},\ Z\> .
\end{align*}
By the above and Definition~\ref{D-TTSSR}, we obtain \eqref{E-01}:
\begin{equation*}
 ({\cal T}_{1}+S_{1}+{\cal T}_{2}+S_{2})(Y,X_{1},X_{2},Z) = -{R^P}(Y,X_{1},X_{2},Z).
\end{equation*}
By the above, since ${\cal T}_{i}$ and $S_{i}$ are tensor fields, also ${R^P}$ is a tensor field.
\hfill$\Box$

\begin{remark}\rm
If $P_{1}$ and $P_{2}$ are orthoprojectors from $TM$ onto complementary orthogonal distributions ${\cal D}_1$ and ${\cal D}_2$, respectively, then $P=\id_{TM}$ and
\begin{eqnarray*}
 && {\cal T}_{1}(Y,X_{1},X_{2},Z)=\<P_{2}(\nabla_{P_{1}X_{1}}{B}_{2})(X_{2}, Y),\ Z\>,\\
 && {\cal T}_{2}(Y,X_{1},X_{2},Z)=\<P_{1}(\nabla_{P_{2}Y}{B}_{1})(Z, X_{1}),\ X_{2}\>,\\
 && S_{1}(Y,X_{1},X_{2},Z)=\<{B}_{2}(X_{2}, {B}_{2}(X_{1},Y)),\ Z\>, \\
 && S_{2}(Y,X_{1},X_{2},Z)=\<{B}_{1}(Z, {B}_{1}(Y,X_{1})),\ X_{2}\>,\\
 && {R^P}(Y,X_{1},X_{2},Z)= {R}(P_2Y, P_1X_{1}, P_1X_{2}, P_2Z).
\end{eqnarray*}
Thus, \eqref{E-01} yields the {Codazzi equation}, see \cite[Lemma~2.25]{rov-m},
\begin{eqnarray}\label{E-Cod}
\nonumber
&& {R}(P_2Y, P_1X_{1}, P_1X_{2}, P_2Z)
+\<(\nabla_{X_1}{B}_{2})(X_2, Y), Z\> +\<(\nabla_{Y}{B}_{1})(Z, X_1), X_2\>\\
&& +\,\<{B}_{2}(X_2,{B}_{2}(X_1,Y)), Z\> +\<{B}_{1}(Z,{B}_{1}(Y,X_1)), X_2\>
=0,
\end{eqnarray}
where $X_i\in{\mathcal D}_{1}$, $Y,Z\in{\mathcal D}_{2}$.
For $B_1=0$, $X_i=X$, \eqref{E-Cod} yields the {Riccati equation}
\[
 (\nabla_{X}{B}_{2})(X, Y) +{B}_{2}(X,{B}_{2}(X,Y)) + {R}(P_2Y, P_1X) P_1X =0.
\]
\end{remark}

\begin{example}\rm
We show the existence of allowed $P=P_1+P_2$ in some cases.
We say that $P_{1},P_{2}\in\mathrm{End}(TM)$ give a \emph{local split}
of $V=U{\times}\bar{U}\subset M$ if the following property~holds:

\smallskip
\textbf{S}$_{1}$: $P_{1}(TV)$ is tangent to $\calf$ and
$P_{2}(TV)$ is tangent to $\bar\calf$, when restricted to $V=U\times\bar{U}$,
where $\calf$ and $\bar\calf$ are simple foliations with leaves $U$ and $\bar{U}$, respectively.

\smallskip\noindent
We say that $P_{1}$ and $P_{2}$ give a $\nabla$-\emph{local split} of
$V=U\times\bar{U}$ as in \textbf{S}$_{1}$, if in addition to \textbf{S}$_{1}$,
the following condition~holds:

\smallskip
\textbf{S}$_{2}$:
The connection $\nabla$ restricts to Levi-Civita connections along the leaves of $\calf$ and $\bar\calf$,
 that is $\nabla_XY$ belongs to $T\calf$ when $X,Y\in T\calf$,
 and $\nabla_XY$ belongs to $T\bar\calf$ when $X,Y\in T\bar\calf$.

\smallskip
We say that $P_{1}$ and $P_{2}$ give $\nabla$-\emph{split} of
$P=P_{1}+P_{2}$ if there is an open cover of
domains
$V=U\times\bar{U}$, where $P_{1}$ and $P_{2}$ give a $\nabla$-local split of $P$.
We conclude with the claim:
 \textit{If $P_{1}$ and $P_{2}$ give a $\nabla$-split of $V\subset M$, then $P=P_{1}+P_{2}$ is allowed for the
 Levi-Civita connection~$\nabla$}.
\end{example}

\section{The modified divergence}
\label{sec:02}

Here, we assume that
 $P\in\mathrm{End}(TM)$ is allowed for the Levi-Civita connection $\nabla$ of metric $g$, and \eqref{E-P1P2} holds.
We extend the divergence formula for vector and tensor fields.
 Recall that the~\textit{divergence} ${\rm div}\,X$ of a vector field $X\in{\cal X}_M$ on
a Riemannian manifold $(M,g)$ is given~by
\begin{equation}\label{E-divf-1}
 d(\iota_{X}\,d\operatorname{vol}) = ({\rm div}X )\,d\operatorname{vol},
\end{equation}
where $d\operatorname{vol}$ is the volume form of $g$ and
$\iota_{\,X}$ is operator of {contraction}.
The divergence of a $(1,k)$-tensor $S$ is a $(0,k)$-tensor
${\rm div}\,S =\operatorname{trace}\left(Y\longrightarrow\nabla_{Y} S\right)$,
that is in coordinates,
\[
 ({\rm div}\,S)_{i_{1},\dots i_{k}} = \nabla_{j}\,S^{j}_{\,\cdot\, i_{1},\dots i_{k}}.
\]

\begin{remark}\rm
Using the equality
\begin{equation}\label{E-d-detg}
 \partial_{i}(\sqrt{\det g})=\sqrt{\det g}\cdot g^{jk}\frac{\partial g_{jk}}{\partial x^{i}}
\end{equation}
and definition of Christoffel symbols, we obtain
\begin{align}\label{E-divf-2}
\nonumber
 {\rm div}\,X  & =X^{i}_{,i}+\frac12\,{g^{ij}}\,g_{ij},_{k}X^{k}\\
 & = X^{i}_{,i}+X^{i}\frac{\partial( \ln\sqrt{\det g}\,)}{\partial x^{i}}
 =\frac{1}{\sqrt{\det g}}\frac{\partial(\sqrt{\det g}\,X^{i})}{\partial x^{i}}.
\end{align}
In coordinates, for a (1,1)-tensor $S$ we have
 $(\nabla_{\partial_{i}}S)(\partial_{j}) = \big( S^{k}_{j,\,i}
 +S_{j}^{l}\Gamma_{ij}^{k} -\Gamma_{ij}^{l}S_{l}^{k}\big)\partial_{k}$,
where $\partial_{i}=\frac{\partial}{\partial x^{i}}$.
Then, using \eqref{E-d-detg}, we get
\begin{eqnarray}\label{MainFormula-A}
\nonumber
 ({\rm div}\,S)_{j} \eq S^{i}_{j,\,i} +S_{j}^{l}\Gamma_{il}^{i}-\Gamma_{ij}^{l}S_{l}^{i}
 = S^{i}_{j,\,i} -\frac12\,S^{ik}\Big(\frac{\partial g_{ik}}{\partial x^{j}}
 -g_{jk}g^{ql}\frac{\partial g_{ql}}{\partial x^{i}}\Big)\\
 \eq \frac{1}{\sqrt{\det g}}\frac{\partial(\sqrt{\det g}\,S_{j}^{i})}{\partial x^{i}}
 -\frac{1}{2}\,S^{ik} \frac{\partial g_{ik}}{\partial x^{j}}.
\end{eqnarray}
\end{remark}

\begin{definition}\label{D-divP}\rm
Given $P\in \mathrm{End}(TM)$, the $P$-\textit{divergence} of a $(1,k)$-tensor $S$ is a $(0,k)$-tensor
\begin{equation*}
 {\rm div}_{P}\,S=\operatorname{trace}(Y\rightarrow P^*\nabla_{P Y}\,S),
\end{equation*}
e.g. for a vector field $X$ on $M$ we get a function
 ${\rm div}_{P}\,X = \operatorname{trace}(Y\rightarrow P^*\nabla_{P Y}\,X)$ on $M$.
\end{definition}

\begin{lemma}
For $P\in\mathrm{End}(TM)$ and any vector field $X$ on $M$, we have in coordinates
\begin{equation}\label{E-div-PX}
 {\rm div}_{P}\,X = (PP^*)_{j}^{i}\,X^{j},_{i} +\frac{1}{2}\,(PP^*)^{ij}\,g_{ij},_{k} X^{k}.
\end{equation}
\end{lemma}

\proof
Given $X\in{\cal X}_M$, the map $Y\longrightarrow P^{\ast}\nabla_{PY}X$ has the local form
\[
 \partial_{i}\longrightarrow P^{\ast}\nabla_{P\partial_{i}}(X^{k}\partial_{k})
 =P_{i}^{l}\big(X^{k}_{,\,l}+X^{s}\Gamma_{ls}^{k}\big)(P^*)_{k}^{j}\,\partial_j,
\]
where $\partial_{i}=\dfrac{\partial}{\partial x^{i}}$ and $(P_{j}^{i})$ -- the components of $P$.
The trace of the above map is
\begin{equation*}
 \operatorname{div}_{P}X =P_{i}^{l}\left(X^{k}_{,\,l}+X^{s}\Gamma_{ls}^{k}\right)(P^*)_{k}^{i}
 =(PP^*)_{k}^{l}\left(X^{k}_{,\,l} +X^{s}\Gamma_{ls}^{k}\right).
\end{equation*}
By the above, using the symmetry of $PP^*$, i.e., $\<PP^{*}(X),Y\>=\<X,PP^{*}(Y)\>$,
and definition of Christoffel symbols $\Gamma_{ik}^{j}$, we~get \eqref{E-div-PX}.
\hfill$\Box$

\begin{proposition}\label{P-03}
Given $P\in \mathrm{End}(TM)$, condition
\begin{equation}\label{E-cond-PP}
 {\rm div}(P P^{*})=0
\end{equation}
is equivalent to the following:
\begin{equation}\label{E-divf-4}
 {\rm div}_{P}\,X={\rm div} (PP^{*}(X)),\quad X\in{\cal X}_{M},
\end{equation}
which means that $({\rm div}_{P}X)\,d\operatorname{vol}$ is an exact form:
\begin{equation}\label{E-divf-5}
 ({\rm div}_{P}\,X)\,d\operatorname{vol} =
 d(\iota_{PP^{*}(X)}\,d\operatorname{vol}).
\end{equation}
Moreover, we have
\begin{equation}\label{E-divf-6}
 {\rm div}_{P}\,X=\<PP^{*},\,\nabla X\>,\quad X\in{\cal X}_{M}.
\end{equation}
\end{proposition}

\proof
From the definition of ${\rm div}_{P}\,X$ and \eqref{E-divf-2}, \eqref{MainFormula-A}, we have for $S=PP^{*}$:
\begin{align*}
 {\rm div}(S(X)) & \overset{\eqref{E-divf-2}}
 = \frac{1}{\sqrt{\det g}} \frac{\partial(\sqrt{\det g}\,S_{j}^{i}X^{j})}{\partial x^{i}}
 = S_{j}^{i}\,X^{j}_{,\,i}
 +\frac{1}{\sqrt{\det g}}\frac{\partial(\sqrt{\det g}\,S_{j}^{i})}{\partial x^{i}}\,X^{j}\\
 &\overset{\eqref{MainFormula-A}}
 =S_{j}^{i}\,X^{j}_{,\,i} +\frac{1}{2}\,S^{ik}\,\frac{\partial g_{ik}}{\partial x^{j}}\,X^{j}
 +({\rm div}\,S)_{j}X^j \overset{\eqref{E-div-PX}}={\rm div}_{P}\,X +{\rm div}(PP^{*})(X),
\end{align*}
thus the first claim follows.
By \eqref{E-divf-1} and \eqref{E-divf-4}, we obtain \eqref{E-divf-5}.
From the above and identity
\begin{equation*}
 {\rm div}(PP^{*}(X))=\<PP^{*},\,\nabla X\> + \<X, {\rm div}(PP^{*})\>,\quad X\in{\cal X}_{M},
\end{equation*}
follows \eqref{E-divf-6}.
\hfill$\Box$

\begin{remark}\rm
Similar to \eqref{E-divf-4} result can be obtained for a $(1,k)$-tensor $S$.
Conditions \eqref{E-cond-PP}--\eqref{E-divf-6} are satisfied for $P=\id_{TM}$.
\end{remark}

\begin{corollary}
 Let \eqref{E-cond-PP} hold. Then the following $($well-known for $P=\id)$ formula is valid:
\begin{equation*}
 {\rm div}_{P}(f\cdot X) = f\cdot{\rm div} (PP^{*}(X)) + (PP^{*}(X))(f),\quad X\in{\cal X}_{M},\ f\in C^\infty(M).
\end{equation*}
\end{corollary}

From Proposition~\ref{P-03} we obtain the following generalization of Stokes theorem,
which for $P=\id_{TM}$ reduces to the classical divergence theorem.

\begin{theorem}\label{T-Stokes-P}
 If \eqref{E-cond-PP} holds on a compact manifold $(M,g)$, then for any $X\in{\cal X}_M$,
\[
 \int_{M} ({\rm div}_{P}\,X)\,d\operatorname{vol} = \int_{\partial M}\<X, PP^{*}(\nu)\>\,d\omega.
\]
\end{theorem}

Next, we modify Stokes' theorem on a complete open Riemannian manifold $(M,g)$.

\begin{proposition}[see \cite{csc2010} for regular case and $P=\operatorname{id}_{TM}$]\label{L-Div-1}
Let $(M,g)$ be a complete open Riemannian manifold endowed with
a vector field $X$ such that ${\rm div}_{P}\,X\ge0$ $($or ${\rm div}_{P}\,X\le0)$,
where $P\in\mathrm{End}(TM)$ such that \eqref{E-cond-PP} and $\|PP^*(X)\|_{g}\in\mathrm{L}^{1}(M,g)$ hold.
Then ${\rm div}_{P}\,X\equiv0$.
\end{proposition}

\proof
Let $\omega$ be the $(n-1)$-form in $M$ given by $\omega= \iota_{PP^*(X)}\,d\operatorname{vol}_{g}$,
i.e., the contraction of the volume form
$d\operatorname{vol}_{g}$ in the direction of a smooth vector field $PP^*(X)$ on $M$.
If $\{e_{1}, \ldots, e_{n}\}$ is an orthonormal frame on an open set
$U \subset M$, with coframe ${\omega_{1}, \ldots, \omega_{n}}$, then
\[
 \iota_{PP^*(X)}\,d\operatorname{vol}_{g} = \sum\nolimits^{n}_{i=1} (-1)^{i-1}
 \<PP^*(X), e_{i}\>\,\omega_{1}\wedge\ldots\wedge\hat\omega_{i} \wedge\ldots\wedge\omega_{n}.
\]
Since the $(n - 1)$-forms $\omega_{1}\wedge\ldots\wedge\hat\omega_{i}
\wedge\ldots\wedge\omega_{n}$ compose an orthonormal frame in $\Omega^{n-1}(M)$, we get
\[
 \|\omega\|_{g}^{2} = \sum\nolimits^{n}_{i=1} \<PP^*(X), e_{i}\>^{2} =\|PP^*(X)\|_{g}^{2}.
\]
By this and conditions, we get $\|\omega\|_{g}\in\mathrm{L}^{1}(M,g)$ and
$d\omega=d(\iota_{PP^*(X)}\,d\operatorname{vol}_{g})=({\rm div}_{P} X)\,d\operatorname{vol}_{g}$, see~\eqref{E-divf-4}.
According to \cite{yau}, there exists a sequence of domains $B_{i}$ on $M$ such that $M=\bigcup_{\,i\ge1} B_{i}$,
$B_{i}\subset B_{i+1}$ and  $\lim_{\,i\to\infty}\int\nolimits_{B_{i}} d\omega= 0$.
Then we obtain
\[
 \int_{B_{i}} ({\rm div}_{P}\,X)\,d\operatorname{vol}_{g} = \int_{B_{i}}
 {\rm div}(PP^*(X))\,d\operatorname{vol}_{g} = \int_{B_{i}} d\omega\to0.
\]
By conditions and Proposition~\ref{P-03} we find that ${\rm div}_{P}\,X = 0$ on $M$.
\hfill$\square$

\begin{example}\label{Ex-Einstein}\rm
Recall that Einstein tensor is divergence free,
 thus it can play a role of $PP^{*}$.
Consider the product $M^{5}=S^{3}\times T^{2}$ and the coordinates $(x,y,z,\,u,v)$,
where $(x,y,z)\in\RR^3$ are stereographic projections from the north pole of $S^3$ and $(u,v)\in[0,2\pi)^{2}$
are the angular coordinates on $T^{2}=S^{1}\times S^{1}$.
Consider the following Einstein metric $g$ on $M^{5}$:
\[
 ds^{2}=\frac{4}{( x^{2}+y^{2}+z^{2}+1)^{2}}\,(dx^{2}+dy^{2}+dz^{2})+(1+\sin^{2}u)(du^{2}+dv^{2}).
\]
The Einstein tensor has diagonal form $E={\rm diag}(-E_{1},-E_{1},-E_{1},-E_{2},-E_{2})$ with
\[
 E_{1} =-\frac{\sin^{2} u(4\cos^{4}u-5\cos^{2}u+10)}{( 1+\sin^{2}u)^{3}},\quad E_{2}=-3.
\]
Thus, $E$ is divergence free and there is a mixed (1,1)-tensor,
$P=\sqrt{-E}$, i.e., $-E=PP^{*}$, which has diagonal form $P={\rm diag}(a_{1},a_{1},a_{1},a_{2},a_{2})$ with
\[
 a_{1}=\frac{\sin u\sqrt{4\cos^{4}u-5\cos^{2}u+10}}{(1+\sin^{2}u)^{3/2}},\quad a_{2}=\sqrt{3}.
\]
The positive endomorphism $P$ is a sum $P=P_{1}+P_{2}$, where $P_{1}$ and $P_{2}$ have diagonal forms
\[
 P_{1}={\rm diag}(a_{1},a_{1},a_{1},0,0),\quad P_{2}={\rm diag}(0,0,0,a_{2},a_{2}).
\]
We claim that $P$ is allowed for the Levi-Civita connection $\nabla$.
Indeed, consider the Christoffel
symbols $\Gamma_{BC}^{A}$, where $A$, $B$, $C\in\{1,\ldots,6\}$ and
$(x^{1},x^{2},x^{3},\,x^{4},x^{5}) = \left( x,y,z,\,u,v\right)$.
Then one can check that $\Gamma_{BC}^{A}=0$,
provided that $\{B,C\}\subset\{1,2,3\}$ and $A\in\{4,5\}$, or $A\in\{1,2,3\}$ and $\{B,C\}\subset\{4,5\}$.
This implies the claim.
\end{example}

\begin{example}\rm
If $P$ is an almost complex structure, then $PP^{*}=\id_{TM}$ and ${\rm div}_P\,X={\rm div} X$.
This simple observation can be developed as follows.

\smallskip
a)~An \textit{almost contact manifold} $(M,\phi,\xi,\eta)$ is an odd-dimensional manifold $M$, which carries a $(1,1)$-tensor field $\phi$, a (Reeb) vector field $\xi$, and a 1-form $\eta$ satisfying, see \cite{b2010},
\[
 \phi^2 = -\id_{TM} +\,\eta\otimes\xi,\quad \eta(\xi) = 1.
\]
One may show that $\phi\,\xi=0$ and $\eta\circ\phi=0$.
We get an \textit{almost contact metric structure},
if there is metric $g=\<\cdot,\cdot\>$ such that
\[
 \<\phi X,\phi Y\>=\<X,Y\>-\eta(X)\eta(Y)
 \ \Longleftrightarrow\
 \phi^*\phi = \id_{TM} -\,\eta\otimes\xi.
\]
Thus, $\phi^*=-\phi$ restricted on $\ker\eta$, and $\phi^*(\xi)=\xi$. Setting $Y=\xi$ we get $\eta(X)=\<X,\xi\>$.
Hence $\<\xi,\xi\>=1$. We have, using $\nabla_{e_i}X=0$,
\begin{eqnarray*}
 {\rm div}(\phi\,\phi^*)(X) \eq \sum\nolimits_i\<(\nabla_{e_i}(\phi\,\phi^*))(X), e_i\> \\
 \eq -\sum\nolimits_i\<(\nabla_{e_i}(\eta\otimes\xi))(X), e_i\>
 =-\sum\nolimits_i\<(\nabla_{e_i}(\eta(X)\xi)), e_i\>\\
 \eq -\sum\nolimits_i [ e_i(\eta(X))\<\xi,e_i\> + ({\rm div}\,\xi)\,\eta(X) ] \\
 \eq -\xi\<\xi,X\> - ({\rm div}\,\xi)\,\eta(X) = -\<\nabla_\xi\,\xi - ({\rm div}\,\xi)\xi, X\>.
\end{eqnarray*}
Note that $\nabla_\xi\,\xi$ is orthogonal to $\xi$.
Thus, the condition ${\rm div}(\phi\,\phi^*) = 0$, see \eqref{E-cond-PP}, holds if and only if $\xi$ is a geodesic
vector field ($\nabla_\xi\,\xi=0$) and the distribution $\ker\phi$ is harmonic (${\rm div}\,\xi=0$).

\smallskip
b)~An $f$-\textit{structure} (due to Yano, 1961) on a manifold $M$ is a non null $(1,1)$-tensor
$f$ on $M$ of constant rank such that $f^3 +f = 0$, which generalizes the almost complex and the almost contact structures.
It is known that $TM$ splits into two complementary subbundles $\widetilde{\cal D}=f(TM)$ and ${\cal D}=\ker f$,
and that the restriction of $f$ to $\widetilde{\cal D}$ determines a complex structure on it.
An interesting case of $f$-structure occurs when ${\cal D}$ is parallelizable for which there
exist global vector fields $\xi_i,\ i\in\{1,\ldots, p\}$, with their dual 1-forms $\eta^i$, satisfying \cite{gy}
\[
 f^2 = -\id_{TM} +\sum\nolimits_i\eta^i\otimes\xi_i,\quad \eta^i(\xi_j)=\delta^i_j.
\]
A Riemannian metric $g=\<\cdot,\cdot\>$ is compatible, if
\begin{equation*}
 \<f(X), f(Y)\> = \<X,Y\> - \sum\nolimits_i \eta^i(X)\,\eta^i(Y)
 \ \Longleftrightarrow\
 f^* f = \id_{TM} -\sum\nolimits_i \eta^i\otimes\xi_i.
\end{equation*}
Thus, $f^*=-f$ restricted on $\bigcap_i\ker\eta_i$, and $f^*(\xi_i)=\xi_i$.
Setting $Y=\xi_j$ we get $\eta^i(X)=\<X,\xi_j\>$. Hence $\<\xi_i,\xi_j\>=\delta_{ij}$.
Similarly to point b), we obtain
\begin{eqnarray*}
 && ({\rm div}\,f f^*)(X) = \sum\nolimits_i\<(\nabla_{e_i} f f^*)(X), e_i\>
 =-\sum\nolimits_{i,j}\<(\nabla_{e_i}(\eta^j\otimes\xi_j))(X), e_i\> \\
 && = -\sum\nolimits_{i,j} e_i(\eta^j(X))\<\xi_j,e_i\> - \sum\nolimits_{j} ({\rm div}\,\xi_j)\,\eta^j(X)
 = -\sum\nolimits_j\<\nabla_{\xi_j}\,\xi_j + ({\rm div}\,\xi_j)\xi_j, X\>.
\end{eqnarray*}
Note that $\sum_j\nabla_{\xi_j}\,\xi_j$ (belongs to $f(TM)$) is the mean curvature vector of $\ker f$.
From ${\rm div}\,\xi_j=0$ we get $\<H,\xi_j\>=0$, where $H$ is the mean curvature vector of $f(TM)$.
Thus, the condition ${\rm div}(f f^*) = 0$, see \eqref{E-cond-PP}, holds if and only if both distributions,
$f(TM)$ and $\ker f$, are harmonic.
\end{example}

\section{The integral formula}
\label{sec:03}

In this section, we assume that $P_i\ (i=1,2)$ are self-adjoint for adapted metric
(with the Levi-Civita connection $\nabla$), see Remark~\ref{Rem-01};
then \eqref{E-P1P2} means the orthogonality of singular distributions $P_i(TM)$
and Definition~\ref{D-TTSSR} for $R^P$ reads as
\[
  {R}^P(Y,X_{1},X_{2},Z) = \<P_{2}(\nabla_{P_{2}Y}P\,\nabla_{P_{1}X_{1}} -\nabla_{P_{1}X_{1}}P\,\nabla_{P_{2}Y}
  -\nabla_{P\,[P_{2}Y,\ P_{1}X_{1}]})P_{1}X_{2},\ Z\> .
\]
Let $\{e_{i}\}$ be a local orthonormal frame in $M$.

\begin{lemma}
Given $P_1,P_2\in\mathrm{End}(TM)$, we have
\begin{subequations}
\begin{eqnarray}\label{E-lem1a}
 &&\hskip-12mm\sum\nolimits_{s,t}{\cal T}_{1}(e_{t},e_{s},e_{s},e_{t})
 =\sum\nolimits_{s,t}\big(\<\nabla_{P_{1}e_{s}}P_{1}e_{s},\ P_{1}\nabla_{P_{2}e_{t}}P_{2}e_{t}\>
 -P_{1}e_{s}\<P_{1}\nabla_{P_{2}e_{t}}P_{2}e_{t},\ P_{1}e_{s}\> \big),\\
\label{E-lem2a}
 &&\hskip-12mm\sum\nolimits_{s,t}{\cal T}_{2}(e_{t},e_{s},e_{s},e_{t})
 =\sum\nolimits_{s,t}\big(P_{2}e_{t}\< \nabla_{P_{1}e_{s}}P_{2}e_{t},\ P_{1}e_{s}\>
 +\< \nabla_{P_{2}e_{t}}P_{2}e_{t},\ P_{2}\nabla_{P_{1}e_{s}}P_{1}e_{s}\>\big),\\
\label{E-lem3a}
 &&\hskip-12mm\sum\nolimits_{s,t} S_{2}\left(e_{t},e_{s},e_{s},e_{t}\right)
 =\sum\nolimits_{s,t}\<P_{2}\nabla_{P_{1}e_{s}}P_{1}e_{t},\ \nabla_{P_{1}e_{t}}P_{1}e_{s}\>,\\
\label{E-lem4a}
 &&\hskip-12mm\sum\nolimits_{s,t}S_{1}\left(e_{t},e_{s},e_{s},e_{t}\right)
 =\sum\nolimits_{s,t}\<P_{1}\nabla_{P_{2}e_{s}}P_{2}e_{t},\ \nabla_{P_{2}e_{t}}P_{2}e_{s}\> .
\end{eqnarray}
\end{subequations}
\end{lemma}

\proof
First we will prove the equality
\begin{eqnarray}
\label{E-lem2}
 \sum\nolimits_{s,t}\left(  \<P_{1}\nabla_{P_{2}e_{s}}P_{2}e_{t}, \nabla_{P_{2}e_{t}}P_{2}e_{s}\> +\<\nabla_{P_{2}\nabla_{P_{2}e_{t}}P_{1}e_{s}}P_{2}e_{t}, P_{1}e_{s}\> \right)  =0.
\end{eqnarray}
Put $P_{1}e_{s}=\Pi_{s}^{u}e_{u},\ P_{2}e_{t}=\tilde{\Pi}_{t}^{u}e_{u}$ and
\begin{align*}
 &  \nabla_{P_{1}e_{s}}P_{2}e_{t}=\omega_{st}^{u}e_{u},\quad\nabla_{P_{2}e_{t}}P_{1}e_{s}
 =\tilde{\omega}_{ts}^{v}e_{v}= \tilde{\Pi}_{t}^{w}{\bar{\omega}}_{ws}^{v}e_{v},\\
 &\nabla_{P_{2}e_{t}}P_{1}e_{s}=\Omega_{ts}^{u}e_{u},\quad\nabla_{P_{2}e_{t}}P_{2}e_{s}
 =\tilde{\Omega}_{ts}^{v}e_{v}.
\end{align*}
Since $P_{1}$ and $P_{2}$ are self-adjoint and $P_{1}P_{2}=P_{2}P_{1}=0$, we
have
 $\Pi_{s}^{u}=\Pi_{u}^{s}$,
 $\tilde{\Pi}_{t}^{u}=\tilde{\Pi}_{u}^{t}$,
 $\Pi_{s}^{u}\tilde{\Pi}_{u}^{t}=0$,
 $\tilde\Pi_{s}^{u}{\Pi}_{u}^{t}=0$.
We obtain for both two terms of \eqref{E-lem2},
\begin{align*}
 & \hskip-5pt A=\sum\nolimits_{s,t}\< P_{1}\nabla_{P_{2}e_{s}}P_{2}e_{t}, \nabla_{P_{2}e_{t}}P_{2}e_{s}\>
 =\sum\nolimits_{s,t,u,v}\<P_{1}(\tilde{\Omega}_{st}^{v}e_{v}), \tilde{\Omega}_{ts}^{u}e_{u}\>
 =\sum\nolimits_{s,t,u,v}\tilde{\Omega}_{st}^{v}\tilde{\Omega}_{ts}^{u}\Pi_{u}^{v},\\
 & B=\sum\nolimits_{s,t}\<\nabla_{P_{2}\nabla_{P_{2}e_{t}}P_{1}e_{s}}P_{2}e_{t}, P_{1}e_{s}\> =\sum\nolimits_{s,t,u}\<\nabla_{P_{2}(\Omega_{ts}^{u}e_{u})}P_{2}e_{t}, P_{1}e_{s}\>\\
 & \quad=\sum\nolimits_{s,t,u,v}\Omega_{ts}^{u}\tilde{\Omega}_{ut}^{v}\Pi_{s}^{v} =\sum\nolimits_{s,t,u,v}\tilde{\Omega}_{st}^{v}{\Omega}_{tu}^{s}\Pi_{u}^{v}.
\end{align*}
It follows that the left hand side of \eqref{E-lem2} vanishes:
\begin{align*}
 & \hskip-7mm A+B=\sum\nolimits_{s,t,u,v}\tilde{\Omega}_{st}^{v}\Pi_{u}^{v}
 \big(\tilde{\Omega}_{ts}^{u} + {\Omega}_{tu}^{s}\big)\\
 &  =\sum\nolimits_{s,t,u}\<\nabla_{P_{2}e_{s}}P_{2}e_{t}, P_{1}e_{u}\>
 \big(\<\nabla_{P_{2}e_{t}}P_{2}e_{s}, e_{u}\> +\<\nabla_{P_{2}e_{t}}P_{1}e_{u}, e_{s}\>\big) \\
 & =\sum\nolimits_{s,t,u,s_{1},u_{1}}\<\nabla_{e_{s_{1}}}P_{2}e_{t},
 e_{u_{1}}\> \big( \<\tilde{\Pi}_{s_{1}}^{s}\nabla_{P_{2}e_{t}}P_{2}e_{s},
 P_{1}e_{u_{1}}\> +\<\Pi_{u_{1}}^{u}\nabla_{P_{2}e_{t}}P_{1}e_{u},P_{2}e_{s_{1}}\> \big)\\
 & =\sum\nolimits_{t,s_{1},u_{1}}\<\nabla_{e_{s_{1}}}P_{2}e_{t}, e_{u_{1}}\>
 \big( \<\nabla_{P_{2}e_{t}}P_{2}^{2}e_{s_{1}}, P_{1}e_{u_{1}}\>
 +\<\nabla_{P_{2}e_{t}}P_{1}^{2}e_{u_{1}}, P_{2}e_{s_{1}}\>\big)\\
 & =\sum\nolimits_{t,s_{1},u_{1}}\< \nabla_{e_{s_{1}}}P_{2}e_{t}, e_{u_{1}}\>
 \big( \<\nabla_{P_{2}e_{t}}P_{2}^{2}e_{s_{1}}, P_{1}e_{u_{1}}\>
 +\underline{\<\nabla_{P_{2}e_{t}}P_{1}e_{u_{1}}, P_{2}^{2}e_{s_{1}}\>}\big) =0,
\end{align*}
since the expression in the last large parenthesis vanishes for any $s,t$.
For the fourth line in above calculation of $A+B$ we used orthogonality of the distributions, e.g.
\begin{eqnarray*}
 \langle\tilde{\Pi}_{s_{1}}^{s}\nabla_{P_{2}e_{t}}P_{2}e_{s},P_{1}e_{u_{1}}\rangle =\langle\tilde{\Pi}_{s_{1}}^{s}\nabla_{P_{2}e_{t}}(\tilde{\Pi}_{s}^{s_{2}}e_{s_{2}}),P_{1}e_{u_{1}}\rangle\\
 =\langle\nabla_{P_{2}e_{t}}(\tilde{\Pi}_{s_{1}}^{s}\tilde{\Pi}_{s}^{s_{2}}e_{s_{2}}),P_{1}e_{u_{1}}\rangle =\langle\nabla_{P_{2}e_{t}}(P_{2}^{2}e_{s_{1}}),P_{1}e_{u_{1}}\rangle
\end{eqnarray*}
The underlined term of line 5 was obtained using equalities \eqref{E-lemma1}:
\begin{equation*}
 \langle\nabla_{P_{2}e_{t}}P_{1}^{2}e_{u_{1}},P_{2}e_{s_{1}}\rangle
 =\langle P_{2}\nabla_{P_{2}e_{t}}P_{1}^{2}e_{u_{1}},e_{s_{1}}\rangle\overset{\eqref{E-lemma1}}
 =\langle P_{2}^{2}\nabla_{P_{2}e_{t}}P_{1}e_{u_{1}},e_{s_{1}}\rangle
 =\langle\nabla_{P_{2}e_{t}}P_{1}e_{u_{1}},P_{2}^{2}e_{s_{1}}\rangle.
\end{equation*}
By Definition~\ref{D-TTSSR}, we have
\begin{align*}
 {\cal T}_{1}(Y,X_{1},X_{2},Z)
 =\<\nabla_{P_{1}X_{1}}P_{2}\nabla_{P_{2}Y}P_{1}X_{2}-\nabla_{P_{2}Y}P_{1}\nabla_{P_{1}X_{1}}P_{1}X_{2}
 -\nabla_{P_{2}\nabla_{P_{1}X_{1}}P_{2}Y}P_{1}X_{2}, P_{2}Z\>,
\end{align*}
where (using the metric property of $\nabla$)
\begin{align*}
 &  \<\nabla_{P_{1}X_{1}}P_{2}\nabla_{P_{2}Y}P_{1}X_{2}, P_{2}Z\>
 = (P_{1}X_{1})\<P_{2}\nabla_{P_{2}Y}P_{1}X_{2}, P_{2}Z\>
 -\<P_{2}\nabla_{P_{2}Y}P_{1}X_{2}, \nabla_{P_{1}X_{1}}P_{2}Z\>,\\
 & \<\nabla_{P_{2}Y}P_{1}\nabla_{P_{1}X_{1}}P_{1}X_{2}, P_{2}Z\>
 =-\<P_{1}\nabla_{P_{1}X_{1}}P_{1}X_{2}, \nabla_{P_{2}Y}P_{2}Z\>.
\end{align*}
By the above and \eqref{E-lem2}, we have \eqref{E-lem1a}:
\begin{align*}
 & \hskip-6pt \sum\nolimits_{s,t}{\cal T}_{1}(e_{t},e_{s},e_{s},e_{t})
 =-\!\sum\nolimits_{s,t} (P_{1}e_{s})\<P_{1}\nabla_{P_{2}e_{t}}P_{2}^{2}e_{t},e_{s}\> {-}\!\sum\nolimits_{s,t}\<P_{2}\nabla_{P_{2}e_{t}}P_{1}e_{s},\nabla_{P_{1}e_{s}}P_{2}e_{t}\>\\
 & +\sum\nolimits_{s,t}\<P_{1}\nabla_{P_{1}e_{s}}P_{1}e_{s}, \nabla_{P_{2}e_{t}}P_{2}e_{t}\>
 -\sum\nolimits_{s,t}\<\nabla_{P_{2}\nabla_{P_{1}e_{s}}P_{2}e_{t}}P_{1}e_{s}, P_{2}e_{t}\>\\
 & =\sum\nolimits_{s,t}\big(
 \<P_{1}\nabla_{P_{1}e_{s}}P_{1}e_{s}, \nabla_{P_{2}e_{t}}P_{2}e_{t}\>
 -(P_{1}e_{s})\<P_{1}\nabla_{P_{2}e_{t}}P_{2}^{2}e_{t},e_{s}\>
 \big)\\
 & =\sum\nolimits_{s,t}\big(\<P_{1}\nabla_{P_{1}e_{s}}P_{1}e_{s}, \nabla_{P_{2}e_{t}}P_{2}e_{t}\>
 -\underline{(P_{1}e_{s})\<P_{1}\nabla_{P_{2}e_{t}}P_{2}e_{t},P_{1}e_{s}\>}\big).
\end{align*}
The underlined term in above calculation was obtained using equality $b^{(1)}_2(e_{t},e_{t})=0$.
 Similarly, using $\< \nabla_{P_{1}e_{s}}P_{2}e_{t}, P_{1}e_{s}\> {+}\< \nabla_{P_{1}e_{s}}P_{1}e_{s}, P_{2}e_{t}\>=0$,
we get~\eqref{E-lem2a}.
By Definition~\ref{D-TTSSR}, we~have
\[
 S_{1}(Y,X_{1},X_{2},Z) =\<\nabla_{P_{2}\nabla_{P_{2}Y}P_{1}X_{1}}P_{1}X_{2},
 P_{2}Z\>,\quad S_{2}( Y,X_{1},X_{2},Z) =\<\nabla_{P_{1}\nabla_{P_{1}X_{1}}P_{2}Y}P_{2}Z, P_{1}X_{2}\>.
\]
Then, using dual for \eqref{E-lem2}, we get \eqref{E-lem3a}:
\begin{equation*}
 \sum\nolimits_{s,t} S_{2}\left(  e_{t},e_{s},e_{s},e_{t}\right)
 =-\sum\nolimits_{s,t}\<\nabla_{P_{1}\nabla_{P_{1}e_{s}}P_{2}e_{t}}P_{1}e_{s},
 P_{2}e_{t}\> =\sum\nolimits_{s,t}\<P_{2}\nabla_{P_{1}e_{s}}P_{1}e_{t},\nabla_{P_{1}e_{t}}P_{1}e_{s}\>.
\end{equation*}
By symmetry in indices, we get \eqref{E-lem4a}.
\hfill$\Box$

\begin{definition}\rm
The \textit{second fundamental forms} $h_i$, the \textit{integrability tensors} $T_i$
and the \textit{mean curvature vectors} $H_{i}={\rm Trace}_{g}\, h_{i}$ of
of singular distributions are defined by
\begin{eqnarray*}
 h_{1}(X,Y) \eq \frac12\,P_{2}( \nabla_{P_{1}X}P_{1}Y +\nabla_{P_{1}Y}P_{1}X), \quad
 h_{2}(X,Y) = \frac12\,P_{1}( \nabla_{P_{2}X}P_{2}Y+\nabla_{P_{2}Y}P_{2}X) ,\\
 {T}_{1}(X,Y) \eq \frac12\,P_{2}\left( \nabla_{P_{1}X}P_{1}Y-\nabla_{P_{1}Y}P_{1}X\right), \quad
 {T}_{2}(X,Y) = \frac12\,P_{1}\left(\nabla_{P_{2}X}P_{2}Y-\nabla_{P_{2}Y}P_{2}X\right),\\
 H_{1} \eq \sum\nolimits_{s}P_{2}\nabla_{P_{1}e_{s}}P_{1}e_{s},\quad
 H_{2}=\sum\nolimits_{s}P_{1}\nabla_{P_{2}e_{s}}P_{2}e_{s}.
\end{eqnarray*}
The definition of $H_i$ is correct because of orthogonality of distributions $P_i(TM)$.
If the second fundamental form vanishes then certain distribution is called
\textit{totally geodesic}, if the integrability tensor vanishes then
the distribution is \textit{integrable},
and if the second fundamental form and integrability tensor vanish simultaneously
then the distribution is called \textit{autoparallel} (for regular case see \cite{BF}).
If the mean curvature vector vanishes then certain distribution is called \textit{harmonic}.
A~distribution ${\cal D}_{1}$ is called \emph{totally umbilical} if there is
$\alpha:M\to\mathbb{N}$ such that
\begin{equation*}
 P_{2}\nabla_{P_{1}X}P_{1} Y = (1/\alpha)\<P_{1}X, P_{1} Y\>\,H_{1}.
\end{equation*}
\end{definition}

Totally umbilical regular distributions appear on twisted products of Riemanni\-an manifolds.
Observe that
\[
 \<H_{1}, X\> =-\mathrm{Trace}(Y\rightarrow B_{1}(Y,X)),\quad
 \<H_{2}, X\>=-\mathrm{Trace}(Y\rightarrow B_{2}(X,Y)).
\]

\begin{definition}\rm
Define the \textit{square of the $P$-norm} of a vector $X\in P_1(TM)\cup P_2(TM)$ by
\begin{equation}\label{E-P-norm}
 \|X\|^2_{P} =\left\{
 \begin{array}{cc}
   \<P_1(X'),X'\> & {\rm if}\ X=P_1(X')\in P_1(TM), \\
   \<P_2(X'),X'\> & {\rm if}\ X=P_2(X')\in P_2(TM).
 \end{array}\right.
\end{equation}
\end{definition}

\begin{remark}\rm
For a general endomorphism $P=P_1+P_2$, the value of $|X|^2_{P}$ is not positive, but we will not use it without its square. We claim that
definition~\eqref{E-P-norm} is correct.
Indeed, if $X=P_{1}(X^{\prime })=P_{1}(X^{\prime\prime})$, then $\langle X,X^{\prime}\rangle =\langle P_{1}X^{\prime\prime},X^{\prime}\rangle
=\langle X^{\prime\prime},P_{1}X^{\prime}\rangle =\langle X^{\prime\prime},X\rangle =\langle X,X^{\prime \prime}\rangle$.
\end{remark}

In~particular, by \eqref{E-P-norm} we have,
\begin{equation}\label{E-Smix-1a}
 \|H_{2}\|_{P}^{2} =\!\sum\nolimits_{s,t}\<P_{1}\nabla_{P_{2} e_{s}} P_{2}e_{s},\ \nabla_{P_{2} e_{t}} P_{2} e_{t}\>,\quad
 \|H_{1}\|_{P}^{2} =\!\sum\nolimits_{s,t}\<P_{2}\nabla_{P_{1} e_{s}} P_{1} e_{s},\ \nabla_{P_{1}e_{t}} P_{1} e_{t}\>,
\end{equation}
which makes sense, since $H_1\in P_2(TM)$ and $H_2\in P_1(TM)$.
Then we define similarly the ``squares of the $P$-norms" of tensors,
\begin{equation*}
 \| h_{1}\|_{P}^{2} =\sum\nolimits_{\,s,t}\|\,{h}_{1}({e}_s,{e}_t)\|_{P}^2,\quad
 \| T_{1}\|_{P}^{2}=\sum\nolimits_{\,s,t}\|\,{T}_{1}({e}_s,{e}_t)\|_{P}^2,\quad {\rm and\ so\ on},
\end{equation*}
which makes sense, since $h_1=P_2 h_1^{\prime}$ and $T_1=P_2 T_1^{\prime}$,  etc.

\begin{lemma} We have
\[
 {\rm div}_{P_{2}}H_{1}={\rm div}_{P}H_{1}+\| H_{1}\|_{P}^{2},\quad
 {\rm div}_{P_{1}}H_{2}={\rm div}_{P}H_{2}+\| H_{2}\|_{P}^{2}.
\]
\end{lemma}

\proof We use Definition~\ref{D-divP},
\begin{align*}
  {\rm div}_{P}\,X =\sum\nolimits_{s}\< P\nabla_{Pe_{s}}X,e_{s}\>,\quad
  {\rm div}_{P_{1}}X =\sum\nolimits_{s}\< P_{1}\nabla_{P_{1}e_{s}}X,e_{s}\>,
\end{align*}
and equality $H_{2}= P_{1}X_{0}$, where $X_{0}=\sum\nolimits_{s}\nabla_{P_{2}e_{s}}P_{2}e_{s}$.
Thus
\begin{eqnarray*}
 &&\hskip-4mm {\rm div}_{P_{1}}H_{2} -{\rm div}_{P}H_{2} = {-}\sum\nolimits_{s}\!\big(\<P_{2}\nabla_{P_{1}e_{s}}P_{1}X_{0},e_{s}\>
 +\<P_{1}\nabla_{P_{2}e_{s}}P_{1}X_{0},e_{s}\> +\<P_{2}\nabla_{P_{2}e_{s}}P_{1}X_{0},e_{s}\>\big)\\
 && =-\sum\nolimits_{s}\<P_{2}\nabla_{P_{2}e_{s}}P_{1}X_{0},e_{s}\>
 =\sum\nolimits_{s}\<P_{1}X_{0},\nabla_{P_{2}e_{s}}P_{2}e_{s}\>=\< P_{1}X_{0},X_{0}\>
 =\| H_{2}\|_{P}^{2},
\end{eqnarray*}
since $\sum\nolimits_{s}\< P_{2}\nabla_{P_{1}e_{s}}P_{1}X_{0},e_{s}\>
=\sum\nolimits_{s}\<\nabla_{P_{1}e_{s}}P_{1}X_{0},P_{2}e_{s}\>=0$.
Indeed, if $P_{1}e_{s}=\sum\nolimits_{u}\Pi_{s}^{u}e_{t}=\Pi_{s}^{u}e_{u}$
and $P_{2}e_{s}=\bar{\Pi}_{s}^{v}e_{v}$, then
$\sum\nolimits_{s}\Pi_{s}^{u}\bar{\Pi}_{s}^{v}=0\ (1\le u,v\le m)$,
since $\{e_{i}\}_{1\le i\le m}$ is an orthonormal frame and $P_{i}$ are self-adjoint.
This completes the proof for $H_{2}$.
The proof for $H_{1}$ is similar.\hfill$\Box$

\smallskip

The mixed scalar curvature, $S_{\rm mix}$, which is an averaged mixed sectional curvature (a~plane, which intersects nontrivially both distributions, is called mixed), is the simplest curvature invariant
of a Rieman\-nian manifold endowed with two
complementary orthogonal distributions, e.g. \cite{rov-m}. The \textit{mixed
scalar curvature of a pair} $(P_{1},P_{2})$ is defined by
\begin{align*}
 S^P_{\rm mix}=\sum\nolimits_{s,t}{R^P}(e_{t},e_{s},e_{s},e_{t}).
\end{align*}
and coincides with $S_{\rm mix}$ for the regular case of an almost product structure.

The above tensors are involved in the formula below, which for regular case belongs to~\cite{Wa01}.

\begin{proposition}\label{P-04}
Given self-adjoint $P_1,P_2\in\mathrm{End}(TM)$, put $P=P_{1}+P_{2}$. Then we have
\begin{equation}\label{E-PW1}
 {\rm div}_{P}( H_{1}+H_{2}) = S^P_{\rm mix}+\| h_{1}\|_{P}^{2} +\| h_{2}\|_{P}^{2}-\| {T}_{1}\|_{P}^{2}
 -\| {T}_{2}\|_{P}^{2}-\| H_{1}\|_{P}^{2}-\| H_{2}\|_{P}^{2}.
\end{equation}
\end{proposition}

\proof We find
\begin{eqnarray}\label{E-Smix-1b}
\nonumber
 & \hskip-3mm {\rm div}_{P}( H_{1}+H_{2}) = {\rm div}_{P_{1}}H_{2}+{\rm div}_{P_{2}}H_{1} -\| H_{2}\|_{P}^{2}-\| H_{1}\|_{P}^{2}\\
 \nonumber
 & =\sum\nolimits_{s,t}(P_{1}e_{t}\<P_{1}\nabla_{P_{2}e_{s}}P_{2}e_{s},P_{1}e_{t}\>
 -\<P_{1}\nabla_{P_{2}e_{s}}P_{2}e_{s},\nabla_{P_{1}e_{t}} P_{1}e_{t}\>)\\
\nonumber
 & +\sum\nolimits_{s,t}(P_{2}e_{t}\<P_{2}\nabla_{P_{1}e_{s}}P_{1}e_{s},P_{2}e_{t}\>
 -\<P_{2}\nabla_{P_{1}e_{s}}P_{1}e_{s}, \nabla_{P_{1}e_{t}}P_{2}e_{t}\>)\\
 & -\sum\nolimits_{s,t}\<P_{1}\nabla_{P_{2}e_{s}}P_{2}e_{s}, \nabla_{P_{2}e_{t}}P_{2}e_{t}\>
 -\sum\nolimits_{s,t}\<P_{2}\nabla_{P_{1}e_{s}}P_{1}e_{s}, \nabla_{P_{1}e_{t}}P_{1}e_{t}\>
\end{eqnarray}
and
\begin{eqnarray*}
 \| h_{1}\|_{P}^{2} -\| {T}_{1}\|_{P}^{2}
 =\sum\nolimits_{s,t}\<P_{2}\nabla_{P_{1}e_{s}}P_{1}e_{t},\nabla_{P_{1}e_{t}}P_{1}e_{s}\>,\\
 \| h_{2}\|_{P}^{2} -\| {T}_{2}\|_{P}^{2}
 =\sum\nolimits_{s,t}\<P_{1}\nabla_{P_{2}e_{s}}P_{2}e_{t},\nabla_{P_{2}e_{t}}P_{2}e_{s}\>.
\end{eqnarray*}
By the above,
\begin{eqnarray}\label{E-Smix-1c}
\nonumber
 && \|h_{1}\|_{P}^{2} +\|h_{2}\|_{P}^{2} -\|{T}_{1}\|_{P}^{2} -\|{T}_{2}\|_{P}^{2}\\
 && =\sum\nolimits_{s,t}\big(\<P_{2}\nabla_{P_{1}e_{s}}P_{1}e_{t}, \nabla_{P_{1}e_{t}}P_{1}e_{s}\>
 + \<P_{1}\nabla_{P_{2}e_{s}}P_{2}e_{t},\nabla_{P_{2}e_{t}}P_{2}e_{s}\>\big).
\end{eqnarray}
Summing \eqref{E-Smix-1a}, \eqref{E-Smix-1b} and \eqref{E-Smix-1c} and eliminating underlined terms, we have
\begin{eqnarray}\label{E-Smix-1}
\nonumber
 -{\rm div}_{P}(H_{1}+H_{2}) +\|h_{1}\|_{P}^{2} +\|h_{2}\|_{P}^{2} -\|{T}_{1}\|_{P}^{2} -\|{T}_{2}\|_{P}^{2}
 -\| H_{1}\|_{P}^{2} -\| H_{2}\|_{P}^{2} \\
\nonumber
 =-\sum\nolimits_{s,t}\big(P_{1}e_{t}\< P_{1}\nabla_{P_{2}e_{s}}P_{2}e_{s},\,P_{1}e_{t}\>
 -\< P_{1}\nabla_{P_{2}e_{s}}P_{2}e_{s},\,\nabla_{P_{1}e_{t}}P_{1}e_{t}\>\big)\\
\nonumber
 -\sum\nolimits_{s,t}\big(P_{2}e_{t}\< P_{2}\nabla_{P_{1}e_{s}}P_{1}e_{s},\,P_{2}e_{t}\>
 -\< P_{2}\nabla_{P_{1}e_{s}}P_{1}e_{s},\,\nabla_{P_{2}e_{t}}P_{2}e_{t}\>\big)\\
\nonumber
 +\sum\nolimits_{s,t}\underline{\< P_{1}\nabla_{P_{2}e_{s}}P_{2}e_{s},\,\nabla_{P_{2}e_{t}}P_{2}e_{t}\>}
 +\sum\nolimits_{s,t}\underline{\underline{\< P_{2}\nabla_{P_{1}e_{s}}P_{1}e_{s},\,\nabla_{P_{1}e_{t}}P_{1}e_{t}\>}}\\
\nonumber
 +\sum\nolimits_{s,t}(\< P_{2}\nabla_{P_{1}e_{s}}P_{1}e_{t},\,\nabla_{P_{1}e_{t}}P_{1}e_{s}\>
 +\< P_{1}\nabla_{P_{2}e_{s}}P_{2}e_{t},\,\nabla_{P_{2}e_{t}}P_{2}e_{s}\>)\\
\nonumber
 -\sum\nolimits_{s,t}\underline{\underline{\< P_{2}\nabla_{P_{1}e_{s}}\,P_{1}e_{s},\,\nabla_{P_{1}e_{t}}\,P_{1}e_{t}\>}}
 -\sum\nolimits_{s,t}\underline{\< P_{1}\nabla_{P_{2}e_{s}}\,P_{2}e_{s},\,\nabla_{P_{2}e_{t}}\,P_{2}e_{t}\>}\\
\nonumber
 =-\sum\nolimits_{s,t}\big(P_{1}e_{t}\< P_{1}\nabla_{P_{2}e_{s}}P_{2}e_{s},\,P_{1}e_{t}\>
 -\< P_{1}\nabla_{P_{2}e_{s}}P_{2}e_{s},\,\nabla_{P_{1}e_{t}}P_{1}e_{t}\>\big)\\
\nonumber
 -\sum\nolimits_{s,t}\big(P_{2}e_{t}\< P_{2}\nabla_{P_{1}e_{s}}P_{1}e_{s},\,P_{2}e_{t}\>
 -\< P_{2}\nabla_{P_{1}e_{s}}P_{1}e_{s},\,\nabla_{P_{2}e_{t}}P_{2}e_{t}\>\big)\\
 +\sum\nolimits_{s,t}\big(\< P_{2}\nabla_{P_{1}e_{s}}P_{1}e_{t},\,\nabla_{P_{1}e_{t}}P_{1}e_{s}\>
 +\< P_{1}\nabla_{P_{2}e_{s}}P_{2}e_{t},\,\nabla_{P_{2}e_{t}}P_{2}e_{s}\>\big).
\end{eqnarray}
Tracing Codazzi equation \eqref{E-01} and using \eqref{E-lem1a}--\eqref{E-lem4a}, we obtain
\begin{eqnarray}\label{E-Smix-2}
\nonumber
 && -S^P_{\rm mix}=\sum\nolimits_{s,t}\left( {\cal T}_{1}+{\cal T}_{2}+S_{1}+S_{2}\right)(e_{t},e_{s},e_{s},e_{t})\\
\nonumber
 && =\sum\nolimits_{s,t}\big(-P_{1}e_{s}\< P_{1}\nabla_{P_{2}e_{t}}P_{2}e_{t},\,P_{1}e_{s}\>
 +\< \nabla_{P_{1}e_{s}}P_{1}e_{s},\,P_{1}\nabla_{P_{2}e_{t}}P_{2}e_{t}\>\\
\nonumber
 && -P_{2}e_{t}\< P_{2}\nabla_{P_{1}e_{s}}P_{1}e_{s},P_{2}e_{t}\>
 +\< \nabla_{P_{2}e_{t}}P_{2}e_{t},\,P_{2}\nabla_{P_{1}e_{s}}P_{1}e_{s}\>\\
 && +\< P_{1}\nabla_{P_{2}e_{s}}P_{2}e_{t},\,\nabla_{P_{2}e_{t}}P_{2}e_{s}\>
 +\< P_{2}\nabla_{P_{1}e_{s}}P_{1}e_{t},\,\nabla_{P_{1}e_{t}}P_{1}e_{s}\>\big),
\end{eqnarray}
Comparing \eqref{E-Smix-1} and \eqref{E-Smix-2},
completes the proof of \eqref{E-PW1}.
\hfill$\Box$

\smallskip
For general $P\in{\rm End}(M)$, the integral of the $P$-divergence of a vector
field over a closed manifold vanishes if we assume \eqref{E-cond-PP}, see Theorem~\ref{T-Stokes-P}.
Thus, under certain assumption for self-adjoint $P$, the integral over the right hand side of \eqref{E-PW1} vanishes.

\begin{theorem}
Given self-adjoint $P_i\in{\rm End}(TM)$ $(i=1,2)$
on a closed Riemannian manifold $(M,g)$, let
\begin{equation}\label{E-cond-PPi}
 {\rm div}(P^2) = 0
\end{equation}
for $P=P_{1}+P_{2}$.
Then the following integral formula holds:
\[
 \int_{M}\big(S^P_{\rm mix} +\| h_{1}\|_{P}^{2} +\| h_{2}\|_{P}^{2} -\| {T}_{1}\|_{P}^{2}
 -\| {T}_{2}\|_{P}^{2} -\| H_{1}\|_{P}^{2} -\| H_{2}\|_{P}^{2}\big)\,{\rm d}\operatorname{vol} =0.
\]
\end{theorem}

\proof
This follows from Propositions~\ref{P-03} and \ref{P-04}.
\hfill$\square$

\smallskip
In the sequel we suppose that $P_i$ are (self-adjoint and) non-negative.

The next results on autoparallel distributions yield splitting of manifolds in regular~case.

\begin{theorem}
Let distributions $P_i(TM)$ be integrable with $H_{1}=0$ on a complete open Riemannian manifold $(M,g)$,
and the leaves $(M^{\prime},g^{\prime})$ of $P_1(TM)$ satisfy
condition $\|H_{2\,|M^{\prime}}\|_{g^{\prime}}\in\mathrm{L}^{1}(M^{\prime},g^{\prime})$,
e.g. $(M^{\prime},g^{\prime})$ are compact, and ${\rm div}(P_1^2) = 0$.
If $S^P_{\rm mix}\ge0$ then $S^P_{\rm mix}\equiv0$ and the distributions are autoparallel.
\end{theorem}

\proof
By conditions, we get
\begin{equation*}
 {\rm div}_{P_{1}}\, H_{2} = S^P_{\rm mix} +\|h_{1}\|_{P}^{2}+\|{h}_{2}\|_{P}^{2}.
\end{equation*}
Using Proposition~\ref{L-Div-1} for each leaf (a complete open manifold), and
since $S^P_{\rm mix}\ge0$ (and $P_i$ are non-negative), we get ${\rm div}_{P_1}\,H_{2}=0$. Thus, ${h}_{i}=0$.
\hfill$\square$

\begin{theorem}
Let distributions $P_i(TM)$ on a complete open Riemannian manifold $(M,g)$ sa\-tisfy $H_{i}=0$ and $T_{i}=0$.
 If $\,S^P_{\rm mix}\ge0$ then $S^P_{\rm mix}\equiv0$ and the distributions are autoparallel.
\end{theorem}

\proof
Under assumptions of our Theorem, we~get
 ${\rm div}_{P} ({H}_{1} +{H}_{2}) = S^P_{\rm mix} +\|{h}_{1}\|_{P}^{2} +\|h_{2}\|_{P}^{2}$.
By Proposition~\ref{L-Div-1} and since $S^P_{\rm mix}\ge0$ (and $P_i$ are non-negative), we get
${\rm div}_{P} ({H}_{1} +{H}_{2})=0$. Thus, ${h}_{1} = h_{2} = 0$.
\hfill$\square$

\smallskip

The next result generalizes \cite[Theorem~4]{step1}.

\begin{theorem}
Let the sets, where the ranks of distributions $P_{1}$ and $P_{2}$ are at least $2$, are dense in
a complete open Riemannian manifold $(M,g)$,
and $\Vert P^2({H}_{1}+{H}_{2})\Vert_{g}\in \mathrm{L}^{1}(M,g)$ for $P=P_{1}+P_{2}$ and \eqref{E-cond-PPi} hold.
Suppose that there exist endomorphisms $Q_{1}$ and $Q_{2}$ such that $Q_{i}^2=P_{i}\ (i=1,2)$,
and the pairs of distributions $(P_{1},Q_{2})$ and $(Q_{1},P_{2})$ are totally umbilical.
If $\,S^P_{\rm mix}\leq 0$ then $S^P_{\rm mix}\equiv0$ 
and the distributions $P_i(TM)$ are autoparallel.
\end{theorem}

\proof By conditions,
\begin{equation}\label{defHH}
 Q_{2}\nabla_{P_{1}X}P_{1}Y =(1/\alpha_{1})\,\<P_{1}X,P_{1}Y\>H_{1,Q_{2}}, \quad
 Q_{1}\nabla_{P_{2}X}P_{2}Y =(1/\alpha_{2})\,\<P_{2}X,P_{2}Y\>H_{2,Q_{1}}.
\end{equation}
We have
\begin{eqnarray*}
 \Vert {h}_{1}\Vert_{P}^{2}-\Vert {T}_{1}\Vert_{P}^{2}
 \eq\sum\nolimits_{s,t}\<P_{2}\nabla_{P_{1}e_{s}}P_{1}e_{t},\nabla_{P_{1}e_{t}}P_{1}e_{s}\>\\
 =\sum\nolimits_{s,t}\<Q_{2}\nabla_{P_{1}e_{s}}P_{1}e_{t},Q_{2}\nabla_{P_{1}e_{t}}P_{1}e_{s}\>
 \eq (\alpha_{1})^{-2}\| H_{1,Q_{2}}\|^{2} \sum\nolimits_{s,t}\<P_{1}e_{s},P_{1}e_{t}\>^{2}.
\end{eqnarray*}
Similarly,
\[
 \|{H}_{1}\|_{P}^{2}=(\alpha_{1})^{-2}\| H_{1,Q_{2}}\|^{2}
 \sum\nolimits_{s,t}\<P_{1}e_{s},P_{1}e_{s}\>\<P_{1}e_{t},P_{1}e_{t}\>.
\]
By the Cauchy-Schwarz inequality (and since $P_i$ are non-negative),
 $\Vert {h}_{1}\Vert_{P}^{2}-\Vert {T}_{1}\Vert_{P}^{2}-\|{H}_{1}\|_{P}^{2}\le 0$.
By symmetry, $\Vert {h}_{2}\Vert_{P}^{2}-\Vert {T}_{2}\Vert_{P}^{2}-\|{H}_{2}\|_{P}^{2}\leq 0$.
 By conditions, from Proposition~\ref{P-04} we~get
 ${\rm div}_{P}\left( H_{1}+H_{2}\right) - S^P_{\rm mix} \le0$.
By this, Proposition~\ref{L-Div-1} and condition $S^P_{\rm mix}\leq 0$,
we get ${\rm div}_{P}({H}_{1}+{H}_{2})=0$ and vanishing of $H_{1,Q_{2}}$ and $H_{2,Q_{1}}$.
Then, using (\ref{defHH}), the conclusion follows.
\hfill $\square $

\end{document}